\documentclass[lettersize,journal]{IEEEtran}
\usepackage{amsmath,amsfonts}
\usepackage{algorithmic}
\usepackage{algorithm}
\usepackage{array}
\usepackage[caption=false,font=normalsize,labelfont=sf,textfont=sf]{subfig}
\usepackage{textcomp}
\usepackage{stfloats}
\usepackage{url}
\usepackage{verbatim}
\usepackage{graphicx}
\usepackage{diagbox}
\usepackage{cite}
\def\BibTeX{{\rm B\kern-.05em{\sc i\kern-.025em b}\kern-.08em
    T\kern-.1667em\lower.7ex\hbox{E}\kern-.125emX}}
\usepackage{balance}
\newtheorem{lem}{Lemma}[section]
\newtheorem{thm}[lem]{Theorem}

\newtheorem{cor}[lem]{Corollary}

\newtheorem{defi}[lem]{Definition}

\begin{document}
\title{Probability model of edge-fault tolerance for regular graphs with respect  to edge connectivity}

\author{Huanshen Jia, Jianguo Qian
\thanks{Manuscript created March, 2025. (Corresponding author: Jianguo Qian.)}
\thanks{Huanshen Jia is with School of Mathematical Sciences, Xiamen University, Xiamen 361005, China, and also with School of Mathematics and Statistics, Qinghai Minzu University, Xining, Qinghai 810007, China (e-mail:jhs1989@163.com)}
\thanks{Jianguo Qian is with School of Mathematical Sciences, Xiamen University, Xiamen 361005, and also with School of Mathematics and Statistics, Qinghai Minzu University, Xining, Qinghai 810007, China (e-mail:jgqian@xmu.edu.cn)}}

\markboth{Journal of \LaTeX\ Class Files,~Vol.~18, No.~9, September~2020}%
{How to Use the IEEEtran \LaTeX \ Templates}

\maketitle

\begin{abstract}
 We consider the probability model of edge-fault tolerance of a network in the sense of connectivity with link faults. Using graph-theoretical notation, we define the edge-fault (EF) and Menger-type edge-fault (MEF) tolerances of a graph as the probabilities that the graph is connected and strongly Menger edge-connected when each edge has a certain failure probability, respectively. We derive an upper bound on the EF tolerance for regular graphs, which reveals an asymptotical behavior when graphs and edge failure probability are large enough. We also perform a simulation experiment on a number of randomly generated regular graphs and some typically well-used graphs. The numerical results show that,  in addition to their well-structured properties for networks, Hypercubes, M\"{o}bius Cubes, Ary-Cubes and Circulant graphs have also higher EF and  MEF tolerance in general. In particular, the M\"{o}bius Cube has both the highest EF and  MEF tolerance among all involved graphs. The numerical results also hint that,  in contrast to  MEF tolerance,  the EF tolerance of regular graphs is not strongly effected by the graph structure.
\end{abstract}


\section{Introduction}
\IEEEPARstart{T}{he} failure of nodes/processors or links is generally inevitable in large scale networks and computing systems. Therefore, evaluating the reliability of multiprocessor systems has become an essential issue. To evaluate the reliability, various types of fault tolerance models have been introduced and played important roles, among which the fault tolerance with respect to connectivity reflects the connectivity properties of a network with node or link failures and, hence, is one of the core measurements \cite{Takabe,Li1,Li,Jia,Nikoletseas,Li2}.

We focus on the connectivity involving link failures. Typically, the topology structure of a network are modelled as a graph in which the vertices represent nodes or processors,   and the edges represent communication links.  In graph theory, many types of connectivity were introduced from both theoretical and practical sides.  One of the frequently used is the traditional edge-connectivity, which is defined as the maximum number of faulty edges that guarantees the connection of the graph. Precisely, a graph is called $k$-edge connected if the removal of any $s$ edges from the graph with $s<k$ does not disconnect the graph, or by Menger's theory, any two vertices are connected by at least $k$ edge-disjoint paths \cite{Menger}. Therefore, for a network with larger edge-connectivity, more efficient routing can be achieved using edge-disjoint paths, providing parallel routing and high fault tolerance, increasing the efficiency of data transmission, and decreasing transmission time \cite{Chen}.

Another type of connectivity is the Menger-type edge-connectivity \cite{Li1,Li2,Qiao,Jia}. A graph is called strongly Menger edge-connected if any two vertices are connected by as many paths as their smaller
degree.  In this sense, a strongly Menger edge-connected network has relatively well-balanced connectivity property with respect to vertex degrees, which not only guarantee the  traditional connectivity, but a certain number of edge-disjoint paths between any two vertices.

Theoretically, the two types of the connectivity above ensure only the maximum number of  faulty edges that guarantees the connection of a graph. However, in practical applications, even if more edges are faulted, a graph may still remain connected or strongly Menger edge-connected. Therefore, in the sense of  fault tolerance, it is then important to quantify the effect of edge faults or, in the sense of probability model, to consider  the probability that a graph is connected or strongly Menger edge-connected when a given number of edges are faulted? Further, if  each edge has a certain failure probability independent of one another, what is the probability that the graph is connected or strongly Menger edge-connected? However, the number of edge cuts in a graph increases exponentially when the scale of the graph increases. Therefore, to answer this question is even more difficult for a large-scale graph in general. Based on the translation lemma between configurations and random regular graphs provided by
Bollob\'{a}s, in \cite{Nikoletseas}  Nikoletseas et al.  considered this question from the random graph point of view to reveal the asymptotic behaviour of regular graphs.

 For a graph $G$ and a real number $p$ with $0\leq p\leq 1$, we define the {\it edge-fault (EF) tolerance} $t_e(G,p)$ as the probability that $G$ is connected when all the vertices are perfectly reliable and  all edges are failed independently with a certain probability $p$. The {\it Menger-type edge-fault (MEF) tolerance} $t_e^M(G,p)$ is defined analogously, that is,  the probability that $G$ is strongly Menger edge-connected when all edges are failed independently with a certain probability $p$.

Due to the wide use of regular graphs in networks, in this paper we consider the EF tolerance and MEF tolerance of regular graphs. We focus on some typically well-used topology structures, including Hypercube, M\"{o}bius Cube, Ary Cube  \cite{Day,Li,Ashir} and Circulant graph \cite{Boesch1,Beivide,Boesch2,Li3}, all of which have attractive structural properties such as simple labeling scheme, higher connectivity, higher fault-tolerance and routing capabilities, lower regular vertex degree, and therefore, has been widely used as many classical real-world multiprocessor systems \cite{Saad,Kung,Cull,Hsieh,Fan1,Vaidya,Day,Li}. In addition to the nice topological properties above, Hypercube has also a simple intuitive structure,  recursive scalability and symmetry  \cite{Saad}, which has made it one of the most popular topology for the design of multicomputer interconnection networks \cite{Kung}. Efforts to improve some of the nice properties of Hypercube have lead to the evolution of hypercube variants such as the twisted $n$-cube, multiply-twisted $n$-cube and M\"{o}bius cube \cite{Cull,Hsieh,Fan1}, all of which could be classified among the Hypercube-Like graphs (HL-graphs) \cite{Saad,Vaidya} or Generalized Hypercubes \cite{Jia}. Among these Hypercube variants, M\"{o}bius Cube has shown some superior static performance compared to Hypercube of the same dimension, e.g., possesses about a half of the diameter, two-thirds of the expected distance, Hamilton-connectivity, embeddability of various architectures \cite{Cull,Cheng,Fan1,Hsieh} and higher fault-tolerant pancyclicity \cite{Hsieh,Yang}.

Because of their higher connectivity, it is reasonably to expect that Hypercube, M\"{o}bius Cube, Ary Cube  and Circulant graph would also have nice performance on EF tolerance and MEF tolerance. We will use  simulation experiment for our probability fault model and give a numerical analysis in an attempt to realize whether these topology structures have indeed higher EF tolerance and MEF tolerance? To this end, we also perform our simulation on a number of randomly generated regular graphs.

 The rest of this paper is organized as follows. In the following section, we give some necessary definitions and related known results. In the third section, we use a combinatorial approach to derive an upper bound on $t_e(G,p)$ for regular graphs, which reveals an asymptotical behavior when the scale of the graph and the edge failure probability are large enough. In the fourth section, we establish an algorithm for our simulation on the values of  $t_e(G,p)$ and  $t_e^M(G,p)$. As an application, we perform our simulation  on Hypercube (Ary Cube), M\"{o}bius Cubes,  Circulant graphs and 50 randomly generated regular graphs of the same scale. According to our numerical results, the fifth section concludes that, in addition to their well-structured properties for network models,  the  Hypercubes, M\"{o}bius Cubes and Circulant graphs, as reasonably expected, have also both higher EF and  MEF tolerance in general. In particular, the M\"{o}bius Cube has both the highest EF tolerance and  MEF tolerance among all involved graphs. Even so, it is a little surprising  that the difference of the EF tolerances among all the involved graphs is very small. In contrast, the thing changes for MEF tolerance. This hints that, in contrast to  MEF tolerance,  the EF tolerance of regular graphs would be not strongly effected by the graph structure. Finally, we also point out that the graph attaining the maximum values of $t_e(G,p)$ and $t_e^M(G,p)$  among all involved graphs, i.e., the  M\"{o}bius Cube, is not symmetric, implying that the symmetry of a graph is not a necessary structural feature that guarantees high fault tolerance of the graph.
 %Finally, we propose three questions concerning the asymptotical behaviors and the existence of the phase transition threshold for the two tolerances.

\section{Preliminaries}
\noindent In the following we use graph theoretical terminologies. For a vertex $v$ of a graph $G$,  we denote by $\deg_G(v)$ or, simply $d(v)$ without confusion,  the degree of $v$ in $G$, i.e., the number of the vertices adjacent to $v$. In particular, if $\deg_G(v)=k$ for every vertex $v$, then we call $G$ $k$-{\it regular}. The  minimum degree of a graph $G$ among all vertices is denoted by $\delta (G)$. For a set $F$ of vertices or edges, we use $G-F$ to denote the graph obtained from $G$ by deleting $F$. For two vertices $v_0$ and $v_k$, a $(v_0, v_k)$-{\it path} of length $k$ is a finite sequence of distinct vertices $v_0v_1\cdots v_{k}$ such that $v_iv_{i+1}$ is an edge for $0\leq i\leq k-1$. For two vertices $u$ and $v$, a set $F$ of edges is an $(u, v)$-{\it edge cut} if $G-F$ has no $(u, v)$-path. The {\it edge-connectivity} $\lambda (G)$ of $G$ is the minimum cardinality of an edge set $F$ such that $G-F$ is disconnected. The {\it vertex connectivity} $\kappa(G)$ of $G$ is defined analogously. A graph $G$ is called {\it vertex symmetric} if it is vertex transitive, that is, for any two vertices, $G$ has an automorphism that maps one vertex to another.

   The following result is the well known Menger's Max-Flow Min-Cut Theorem (edge version), which will play a key role in the following discussion.

     \begin{thm}\label{thm 2.1}{\rm({\bf Max-Flow Min-Cut Theorem} {\rm \cite{Menger}})}
		For any two  vertices $u$ and $v$ in a graph $G$, the minimum size of an $(u, v)$-edge cut in $G$ equals the maximum number of edge-disjoint $(u, v)$-paths.
	\end{thm}

We now give the formal definitions related to Menger-type edge-connectivity.

     \begin{defi}\label{defi 2.2}{\rm \cite{Li1,Qiao}}
		A connected graph $G$ is called strongly Menger edge-connected if, for any two distinct vertices $u$ and $v$ in $G$, there are $\min\{\deg_G(u),\deg_G(v)\}$-edge disjoint $(u,v)$-paths between $u$ and $v$.
	\end{defi}

	\begin{defi}\label{defi 2.3}{\rm \cite{Li1}}
	For connected graph $G$ and a positive integer $f$, $G$ is $f$-edge-fault-tolerant strongly Menger edge-connected or simply, $f$-strongly Menger edge-connected, if $G-F$ is strongly Menger edge-connected for any edge set $F$ with $|F|\leq f$.
		
	\end{defi}

In the remaining of the section we introduce the definitions of Hypercube, M\"{o}bius cube, Ary Cube and Circulant graph and their elementary properties involving the edge-connectivity.

     \begin{defi}\label{defi 3.1}{\rm ({\bf Hypercube} {\rm \cite{Saad}})}
    For a positive integer $n$, an $n$-hypercube $Q_n$ is an $n$-regular graph consisting of $2^n$ vertices labeled by $n$-bits binary strings $x_1x_2\cdots x_n\in\{0,1\}^n$, such that there is an edge between any two vertices if and only if their binary strings differ by one and only one bit.
     \end{defi}

    \begin{lem}\label{lem 3.2.} {\rm\cite{Armstrong}}
		 $\kappa(Q_n)=\lambda(Q_n)= n$ for any $n \geq 1$.	
	\end{lem}

     \begin{thm}\label{thm 3.3.} {\rm \cite{Qiao}}
		$Q_n$ is $f$-strongly Menger edge connected if $f\leq 2n-4$ and $n\geq 4$.
	\end{thm}

\begin{defi}\label{defi 2.7.}{\rm({\bf M\"{o}bius Cube} \cite{Cull,Hsieh})}
An  $n$-dimensional M\"{o}bius cube $M_n$ is an $n$-regular graph with $2^n$ vertices, each vertex in which is labelled by an $n$-bits binary string $x_1x_2\cdots x_n\in\{0,1\}^n$ such that the vertex $x_1x_2\cdots x_n$ is adjacent to a vertex $Y$ if and only if $Y$ satisfies one of the following conditions:\\
\begin{enumerate}
\item{$Y=x_{1} \cdots x_{i-1}\overline{x_{i}} \cdots x_n$ if $x_{i-1}=0$ for some $i\in\{1,2,\cdots,n\}$, or} \\
\item{$Y=x_1 \cdots x_{i-1}\overline{x_{i}\cdots x_n}$ if $x_{i-1}=1$ for some $i\in\{1,2,\cdots,n\}$, where $\overline{x_{i}}$ is the complement of $x_i$}.
\end{enumerate}
\end{defi}

We note that, for $i=1$, since the bit $x_{i-1}=x_0$ is not defined, we can assume that $x_0$ is either $0$ or $1$. This means that the rule in the above definition generates two graphs due to the two values of $x_0$. We call the generated graph the {\it 0-M\"{o}bius cube} if $x_0=0$, denoted by 0-$M_n$, and {\it 1-M\"{o}bius cube}  if $x_0=1$, denoted by 1-$M_n$.

\begin{lem}\label{lem 2.8.} {\rm\cite{Vaidya}}
		 $\kappa(M_n)=\lambda(M_n)= n$ for any $n \geq 1$.	
	\end{lem}

    \begin{thm}\label{thm 2.9.}{\rm \cite{Li2}}
		$M_n$ is $f$-strongly Menger edge-connected if $f\leq n-2$ and $n\geq 3$.
	\end{thm}

  \begin{defi}\label{ary}{\rm ({\bf Ary Cube} {\rm \cite{Ashir,Li}})}
  The $k$-ary $n$-cube $Q_n^k$ $(k\geq 2$ and $n\geq 1)$ is a regular graph consisting of $k^n$ vertices labeled by $n$-bits strings $x_1x_2\cdots x_n\in\{0,1,\cdots,k-1\}^n$, such that two vertices are adjacent if and only if they agree on all bits except one, and on that bit they differ by 1 modulo $k$.
\end{defi}

We note that the $2$-ary $n$-cube $Q_n^2$ is exactly the Hypercube $Q_n$ and, hence, is $n$-regular. In general, if $k\geq 3$ then $Q_n^k$ is $2n$-regular. Therefore, the notion of Ary Cube is an extension of Hypercube.
    \begin{lem}\label{lem 2.6.} {\rm\cite{Day,Li}}
	$\kappa(Q_n^k)=\lambda(Q_n^k)= 2n$ for any $n \geq 1$ and $k\geq 3$.	
\end{lem}
\begin{defi}\label{defi 2.10.}{\rm({\bf Circulant Graph} \cite{Boesch2,Li3})}
For positive integers $n_1,n_2,\ldots ,n_k$, the Circulant graph $C_p(n_1,n_2,\ldots ,n_k)$ is a $2k$-regular graph consisting of $p$ vertices labelled by $1,2,\cdots,p$ such that the vertex $i$ is adjacent to $j$ if and only if $j= i\pm n_s({\rm mod}\ p)$ for some $s\in\{1,2,\cdots,k\}$.
\end{defi}

 It is well known that a Circulant graph $C_p(n_1,n_2,\ldots ,n_k)$ is connected if and only if $gcd(p,n_1,n_2,\ldots ,n_k)=1$ \cite{Boesch2}. It is also known that every connected vertex symmetric graph $G$ has edge connectivity $\delta(G)$, and so does every connected Circulant graph.

    \begin{lem}\label{lem 2.11.} {\rm\cite{Boesch2}}
	If a Circulant graph $G$ is connected, then $\kappa(G)=\lambda(G)=\delta(G)$.	
	\end{lem}

\section{An upper bound on EF tolerance of regular graphs}
\noindent The {\it independence number} $i(G)$ of a graph $G$ is the maximum number of pairwise non-adjacent vertices in $G$.
  \begin{thm}\label{lp} For any $d$-regular graph $G$ with independence number $i(G)$, we have
  \begin{equation}\label{lpe}
  t_e(G,p)\leq (1-p^d)^{i(G)}.
  \end{equation}
   \end{thm}
  {\it Proof}. Let $W$ be an independent set of $G$ with independence number $i=i(G)$. For any $v\in W$, if all the edges incident with $v$ are removed from $G$, then the remaining graph is not connected. Let P$(v)$ be the property that  all the edges incident with $v$ are removed from $G$. Thus, if  the remaining graph is connected, then for  any $v\in W$, the property P$(v)$ cannot be hold. We note that the  probability that the edges incident with a certain vertex in $W$ are faulted equals $p^d$. Further, for any two vertices $u,v\in W$, since $W$ is independent, the set of the edges incident with $u$ is disjoint with that incident with $v$. So by the Inclusion-Exclusion principle,
 $$\begin{array}{ccl}
t_e(G,p)&\leq&\sum\limits_{U\subseteq W}(-1)^{|U|}{i\choose |U|}p^{d|U|}\\
&=&{i\choose 0}p^{0}-{i\choose 1}p^{d}+{i\choose 2}p^{2d}-\cdots+(-1)^{i} {i\choose i}p^{id}\\
&=& (1-p^d)^i.
\end{array}
$$
This completes our proof.

We notice that the removal of the edges incident with any certain vertex of a graph will result in a disconnected graph. However, in the proof of Theorem \ref{lp} we consider only the vertices in $W$. This means that the equality in (\ref{lpe}) does not hold in general, with only one exception, i.e., the complete graph with two vertices.
 \begin{cor}\label{lpg} Let  $G_1,G_2,\ldots$ be a infinite sequence of regular graphs and the regular degree $d(G_k)\rightarrow+\infty$ as $k\rightarrow+\infty$.  If $i(G_k)\sim c\alpha^{d(G_k)}$ where $\alpha>1$ is a constant and $p>1/\alpha$, then
 \begin{equation}\label{lpg}
 t_e(G_k,p)\rightarrow 0\ \ {\rm as}\ \ k\rightarrow+\infty.
  \end{equation}
 \end{cor}
 {\it Proof}. Write $d(G_k)=d$.
 {$$(1-p^d)^{i(G_k)}\sim(1-p^d)^{c\alpha^{d}}=\left(1-\frac{1}{\left(1/p\right)^d}\right)^{(1/p)^dp^dc\alpha^{d}}$$
$$=\left(1-\frac{1}{\left(1/p\right)^d}\right)^{(1/p)^dc(\alpha p)^d}
=({\rm e}^{-1})^{c(\alpha p)^d}$$
$$=\left\{\begin{array}{cc}
1,&{\rm if}\ p<1/\alpha,\\
{\rm e}^{-c},&{\rm if}\ p=1/\alpha,\\
0,&{\rm if}\ p>1/\alpha.
\end{array}
\right.
\ \ {\rm as}\ \ k\rightarrow +\infty.$$
The corollary follows directly by (\ref{lpe}), which completes our proof.

Notice that $Q_n$ is $n$-regular and $i(Q_n)=1/2\times 2^n$ since $Q_n$ is bipartite. So by (\ref{lpg}), if $p>1/2$ then  $t_e(Q_n,p)\rightarrow 0$, as $n\rightarrow+\infty$. On the other hand, for a  given $p$ and a random $d$-regular graph $G$ with $N$ vertices, it was proved that, if $d\geq \beta\ln N$, then $t_e(G,p)\rightarrow 1$ as $N\rightarrow+\infty$, where $\beta$ is an appropriate constant \cite{Nikoletseas}. Our result for $Q_n$ implies that $\beta>1/\ln 2$, or almost all $n$-regular graphs with $2^n$ vertices have better  EF tolerance than $Q_n$ does. Further, for vertex-fault tolerance, an asymptotic analysis of Takabe et al.\cite{Takabe} predicts that $1/2$ is the phase transition threshold for almost all random $n$ regular graphs $G$ with $2^n$ vertices, that is, the probability that $G$ is connected tends to 1 and 0 as $n\rightarrow +\infty$ when each vertex has failure probability smaller and greater than $1/2$, respectively. For EF tolerance $t_e(G,p)$, we propose the following question:
\vspace{2mm}\\
{\bf Question 1}. Does there exist the phase transition threshold $p$ of $t_e(G,p)$ for almost all random $n$-regular graphs $G$ with $2^n$ vertices? In particular, if  $p\leq1/2$, what is the value of  $t_e(Q_n,p)$ as $n\rightarrow+\infty$?

\section{Simulation and numerical results}
\noindent In this section, we perform our simulation on 4-regular graphs with 16 vertices, that  is, the Hypercube $Q_4$, the 4-Ary 2-Cube $Q_2^4$,  M\"{o}bius Cubes 0-$M_4$ and 1-$M_4$, and the Circulant Graphs $C_{16}(1,i)$ ($i\in\{2,3,\ldots,7\}$). Moreover, to compare with the above graphs, we also perform our simulation on 50 4-regular graphs with 16 vertices that are generated randomly (by Matlab rand function). Recall that the $2$-ary $4$-cube $Q_4^2$ is exactly the hypercube $Q_4$.

It is clear that every 4-regular graphs with 16 vertices has 32 edges. In our simulation, for each graph $G$ and each integer $f$  with $3\leq f\leq32$,  we choose 2000 graphs, namely the {\it faulty graphs}, obtained from $G$ in which $f$ edges, called the {\it faulty edge set}, are faulted randomly (are deleted randomly). When $f\leq 2$, we choose all the possible edge sets as faulty edge sets.

Before giving our algorithm, we recall some results in graph theory, which are necessary in the algorithm. By Max-Flow Min-Cut Theorem, a graph $G$ is not strongly Menger edge-connected if and only if $G$ has an edge cut $F$ that separates two vertices $u$ and $v$ such that $|F|<\min\{d(u),d(v)\}$. Since all the faulty graphs are obtained from 4-regular graphs, the vertex degrees in which are not greater than 4, meaning that, if such edge cut $F$ exists, it must consist of at most three edges. Further, two vertices $u$ and $v$ are connected if and only if the element $(u,v)$ (the element at the $u$-th row and $v$-th column) in the matrix $A^{n-1}+A^{n-2}$ is not 0, where $A$ is the adjacency matrix of $G$ and $n$ the number of the vertices in $G$. Finally, we recall  that if a graph $G$ with 16 vertices has less than 15 edges, then $G$ is not connected and, hence, not strongly Menger edge-connected.

For an integer $f$ and a graph $G$ that we mentioned above, we denote by $\lambda_f(G)$ and $\lambda^M_f(G)$, or simply $\lambda_f$ and $\lambda_f^M$, the numbers of connected graphs and strongly Menger edge-connected among all the 2000 faulty graphs of $G$ with $f$  faulty edges, respectively. Further, we denote by ${\rm p}_f(G)$ and ${\rm p}_f^M(G)$, or simply ${\rm p}_f$ and ${\rm p}_f^M$, the  probabilities that the 2000 faulty graphs with $f$ faulty edges are connected and strongly Menger edge-connected, respectively. It is clear that ${\rm p}_f(G)=\lambda_f(G)/2000$ and  ${\rm p}_f^M(G)=\lambda_f^M(G)/2000$.

   We now give our algorithm for the values of $\lambda_f(G)$ and $\lambda^M_f(G)$ for each $f$ with $1\leq f\leq 17$ as Algorithm 1.

We perform the algorithm on all the graphs we mentioned above and list some numerical results in Appendix, in which  $R_4$ is a randomly generated graph as illustrated in Fig 1.

\begin{figure}[!t]
\centering
\includegraphics[width=1.6in]{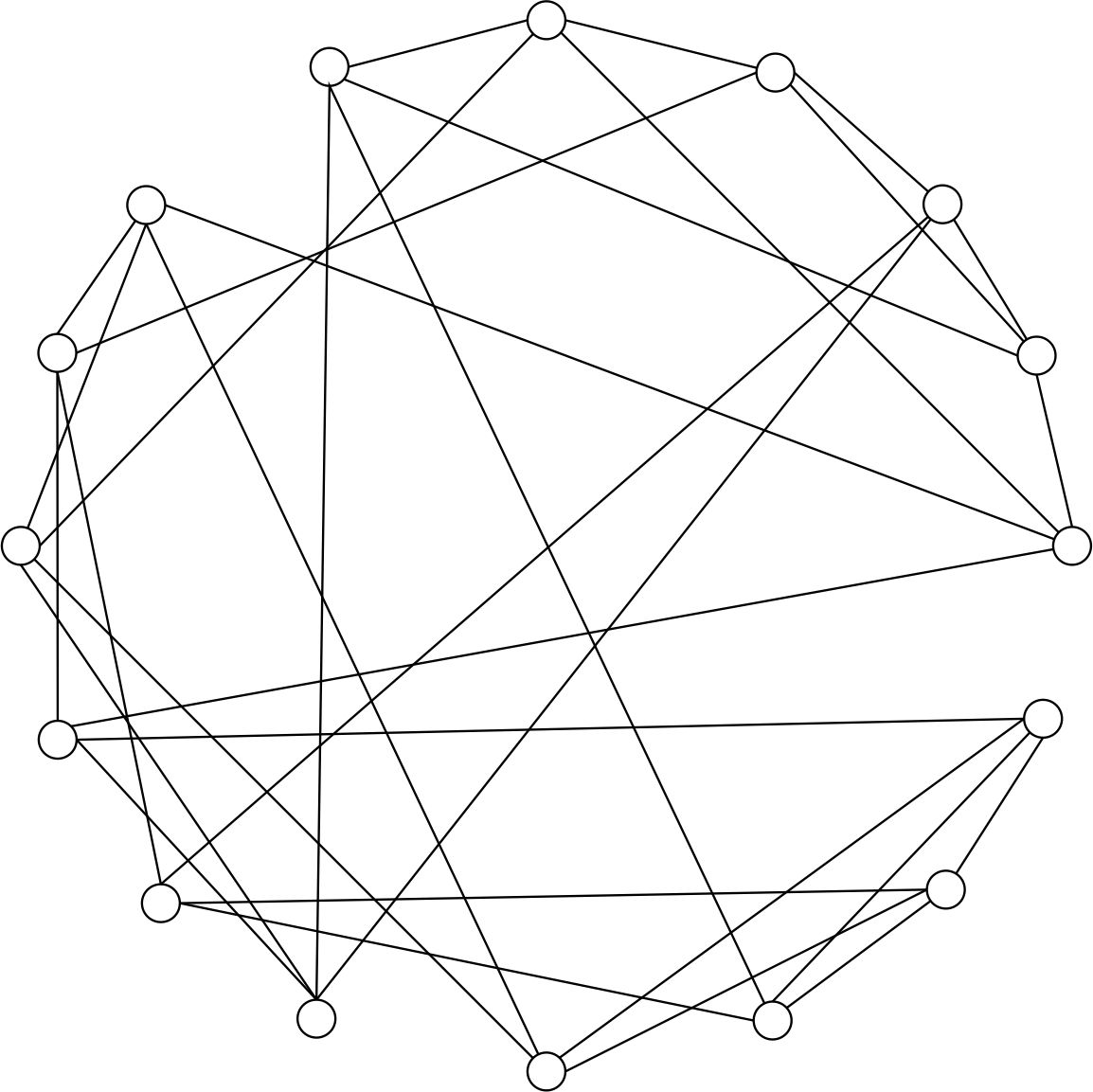}
\caption{A randomly generated 4-regular graph with 16 vertices.}
\label{fig1}
\end{figure}
Intuitively, by Table 1-3, Table 6 and Table 8 in Appendix, we illustrate the probabilities ${\rm p}_f(G)$ and ${\rm p}^M_f(G)$ with $f=0,1,\ldots,18$ as in Fig 2, in which ${\rm p}_f(\overline{H}_4)$ represents the average of  ${\rm p}_f(G)$ among $Q_4$, 0-$M_4$ and 1-$M_4$; ${\rm p}_f(\overline{C}_4)$ the average among all the Circulant graphs $C_{16}(1,i),i=2,3,\cdots,7$; and ${\rm p}_f(\overline{R}_4)$ the average among all the 50 randomly generated graphs. The notation ${\rm p}^M_f(A)$ where $A\in\{\overline{H}_4.\overline{C}_4,\overline{R}_4\}$ is defined analogously to be the average of ${\rm p}^M_f(G)$.
\begin{algorithm}[H]
\caption{The values of $\lambda_f(G)$ and $\lambda^M_f(G)$ for each $f$ with $1\leq f\leq 17$}\label{alg:alg1}
\begin{algorithmic}
\STATE
\STATE \hspace{0.5cm}{\bf Input}: $A$ (the adjacency matrix of $G$);\\
\STATE \hspace{0.5cm}{\bf for} $1\leq f\leq 2$, $F_f=\{F_{f1},F_{f2},\cdots,F_{f{32\choose f}}\}$\\
\STATE \hspace{0.5cm}(the class of all the ${32\choose f}$ faulty edge sets with $f$ edges, represented as matrices);\\
\STATE \hspace{0.5cm}{\bf for} $3\leq f\leq 17$, $F_f=\{F_{f1},F_{f2},\cdots,F_{f2000}\}$\\
\STATE \hspace{0.5cm}(the class of 2000 faulty edge sets with $f$ edges generated randomly from $G$, represented as matrices);\\
\STATE \hspace{0.5cm}{\bf for} $f=1:17$ \\
{\bf Initialization}: $\lambda_{f}=0,\lambda^M_{f}=0$\\
\STATE \hspace{0.5cm}{\bf for} $i=1:2000$\\
\STATE \hspace{0.5cm}{\bf if} $(A-F_{fi})^{14}+(A-F_{fi})^{15}$ has no 0 elements\\
\STATE \hspace{0.5cm}{\bf then} $\lambda_{f}=\lambda_{f}+1$ and let $E_{fi}=:\{E_{fi}^{1},\ldots,E_{fi}^t\}$\\
\STATE \hspace{0.5cm}(the  class of all the edge subsets of $A-F_{fi}$ with at most 3 edges)\\
\STATE \hspace{0.5cm}{\bf for}  $j=1:t$\\
\STATE \hspace{0.5cm}{\bf if}  $(A-F_{fi}-E_{fi}^{j})^{14}+(A-F_{fi}-E_{fi}^{j})^{15}$ has a 0 element $(s,t)$ satisfies\\ \hspace*{24mm}$\min\{d(s),d(t)\}>|E_{fi}^{j}|$ \\
\STATE \hspace{0.5cm}{\bf break}, let $i=i+1$\\
\STATE \hspace{0.5cm}{\bf else} $j=j+1$\\
\STATE \hspace{0.5cm}{\bf if} for all $j$, $(A-F_{fi}-E_{fi}^{j})^{14}+(A-F_{fi}-E_{fi}^{j})^{15}$ has no 0 element $(s,t)$ satisfies\\ \hspace*{19mm}$\min\{d(s),d(t)\}>|E_{fi}^{j}|$\\
\STATE \hspace{0.5cm}{\bf then} $\lambda^M_{f}=\lambda^M_{f}+1$\\
\STATE \hspace{0.5cm}{\bf else} $i=i+1$\\
\STATE \hspace{0.5cm}{\bf output $\lambda_{f}$}\\
\STATE \hspace{0.5cm}{\bf output $\lambda^M_{f}$}\\
{\bf end}
\end{algorithmic}
\label{alg1}
\end{algorithm}

\begin{figure}[!t]
\centering
\includegraphics[width=3.5in]{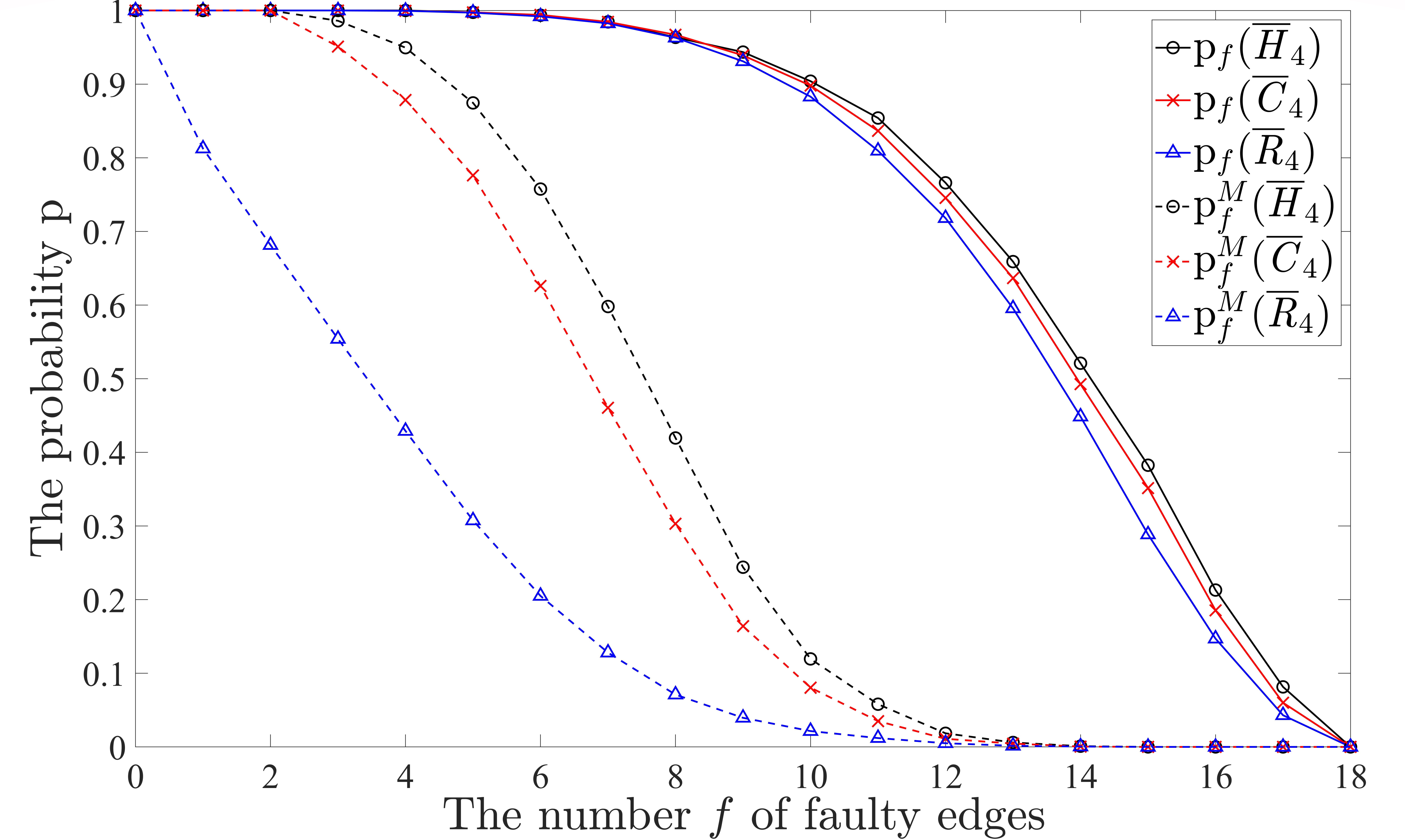}
\caption{The average probabilities ${\rm p}_f(A)$ (solid lines) and ${\rm p}^M_f(A)$ (dotted lines) with $f=0,1,\ldots,18$.}
\label{fig2}
\end{figure}

By the numerical results, we can now consider the two fault tolerances $t_e(G,p)$ and $t_e^M(G,p)$ in terms of $f,{\rm p}_f(G)$ and ${\rm p}^M_f(G)$. Since the failure probability $p$ of each edge is independent from one another, it is clear that
$$t_e(G,p)=\sum\limits_{f=0}^{32}{32\choose f}{\rm p}_f(G)p^f(1-p)^{(32-f)}$$
and
$$t_e^M(G,p)=\sum\limits_{f=0}^{32}{32\choose f}{\rm p}^M_f(G)p^f(1-p)^{(32-f)}.$$

By Table 1-8 in Appendix, we calculate the values of $t_e(G,p)$ and $t_e^M(G,p)$ with $p=0.1,0.2,\ldots,0.9$ for all involved graphs. The values show that $C_{16}(1,4)$ and $C_{16}(1,6)$ attain the maximum values of $t_e(G,p)$ and  $t_e^M(G,p)$ (here, by `the  maximum value' we mean the maximum average value over all $p=0.1,0.2,\ldots,0.9$) among all Circulant graphs, respectively, $R_4$  attains the maximum values of both $t_e(G,p)$ and  $t_e^M(G,p)$ among all randomly generated graphs, and 1-$M_4$  attains the maximum values of both $t_e(G,p)$ and  $t_e^M(G,p)$ among all involved graphs. We list some calculated values in TABLE \ref{tab1} and TABLE \ref{tab2}, in which\\
\begin{itemize}
\item{$\overline{C}^C_4$: the average value of  $t_e(G,p)$ among  all the Circulant graphs $C_{16}(1,i),i=2,3,\ldots,7$;}\\
\item{$\overline{C}^M_4$: the average value of  $t^M_e(G,p)$ among  all the Circulant graphs $C_{16}(1,i),i=2,3,\ldots,7$;}\\
\item{$\overline{R}^C_4$: the average value of  $t^C_e(G,p)$ among  all 50 randomly generated graphs;}\\
\item{$\overline{R}^M_4$: the average value of  $t^M_e(G,p)$  among  all 50 randomly generated graphs.}
\end{itemize}
\begin{table*}
\begin{center}
\caption{The numerical values of $t_e(G,p)$.}
\label{tab1}
{\begin{tabular}{|c|c|c|c|c|c|c|c|c|c|}
		\hline
\diagbox{$p$}{$G$}&$Q_4$   &  0-$M_4$ &  1-$M_4$   &$C_{16}(1,4)$   &$\overline{C}^C_4$  &$R_4$ &$\overline{R}^C_4$\\
		\hline
$0.1$       &0.99838    & 0.99831   & 0.99843    &0.99845    & 0.99836         &0.99859    &0.99807     \\
	\hline
$0.2$       &0.97286   & 0.97306   & 0.97391    &0.97252    & 0.97235         &0.97375     &0.96847     \\
	\hline
$0.3$       &0.86524    &0.86539   & 0.86950    &0.86363    & 0.85813         &0.86398     &0.84196     \\
	\hline
$0.4$       &0.61474	  &0.61486	  &0.62386	  &0.61709	   &0.59927	         &0.61032     &0.56874 \\
	\hline
$0.5$       &0.28818	  &0.28690	  &0.29630	  &0.29182	   &0.27305	          &0.28262	   &0.24643     \\
	\hline
$0.6$       &0.06898    & 0.06799   & 0.07180    &0.07008    & 0.06280         &0.06641     &0.05352     \\
	\hline
$0.7$       &0.00564	    &0.00553	&0.00597	&0.00574	&0.00492	     &0.00534	   &0.00397     \\
	\hline
$0.8$       &0.00007	    &0.00007	&0.00008	&0.00007	&0.00006	     &0.00006	    &0.00004     \\
	\hline
$0.9$       &0            & 0         & 0         &0          & 0              &0             &0     \\
	\hline
{\rm average}&0.42379    & 0.42357   & 0.42665    &0.42438    & 0.41877          &0.42234      &0.40902     \\
	\hline
\end{tabular} }
\end{center}
\end{table*}
\begin{table*}\caption{The numerical values of $t_e^M(G,p)$.}
\begin{center}
\label{tab2}
{\begin{tabular}{|c|c|c|c|c|c|c|c|c|c|}
		\hline
\diagbox{$p$}{$G$}&$Q_4$ &0-$M_4$ &  1-$M_4$  &$C_{16}(1,6)$   & $\overline{C}^M_4$ &$R_4$  &$\overline{R}^M_4$ \\
\hline
$0.1$   &0.94044  &0.94110  &0.94193       &0.94005     &0.89584    &0.89608    &0.53879  \\
\hline
$0.2$   &0.63586	&0.64877	&0.65245	&0.65064	&0.55078	&0.52865	&0.22289  \\
\hline
$0.3$   &0.26013	&0.27466	&0.28095	&0.27666	&0.20893	&0.19282	&0.06487  \\
\hline
$0.4$   &0.06040	&0.06674	&0.06956	&0.06696	&0.04642	&0.04347	&0.01320  \\
\hline
$0.5$   &0.00745	&0.00868	&0.00903	&0.00857	&0.00561	&0.00567	&0.00172  \\
\hline
$0.6$   &0.00041   &0.00051     &0.00052    &0.00048    &0.00031    &0.00035    &0.00011  \\
\hline
$0.7$   &0.00001	&0.00001	&0.00001	&0.00001	&0.00001	  &0.00001	     &0  \\
\hline
$0.8$   &0         &0         &0            &0          &0            &0         &0  \\
\hline
$0.9$   &0         &0         &0            &0          &0            &0         &0  \\
\hline
{\rm average}&0.21163 &0.21561&0.21716      &0.21593    &0.18977     &0.18523   &0.09351  \\
\hline
\end{tabular} }
\end{center}
\end{table*}

By TABLE \ref{tab1} and TABLE \ref{tab2}, we illustrate the values of $t_e(G,p)$ and $t_e^M(G,p)$ for $G\in\{Q_4,$1-$M_4,C_{16}(1,4),$ $C_{16}(1,6),R_4,\overline{R}_4^C,\overline{R}_4^M\}$ and $p=0,0.1,\ldots,0.8$ as in Fig 3.

\begin{figure}[!t]
\centering
\includegraphics[width=3.55in]{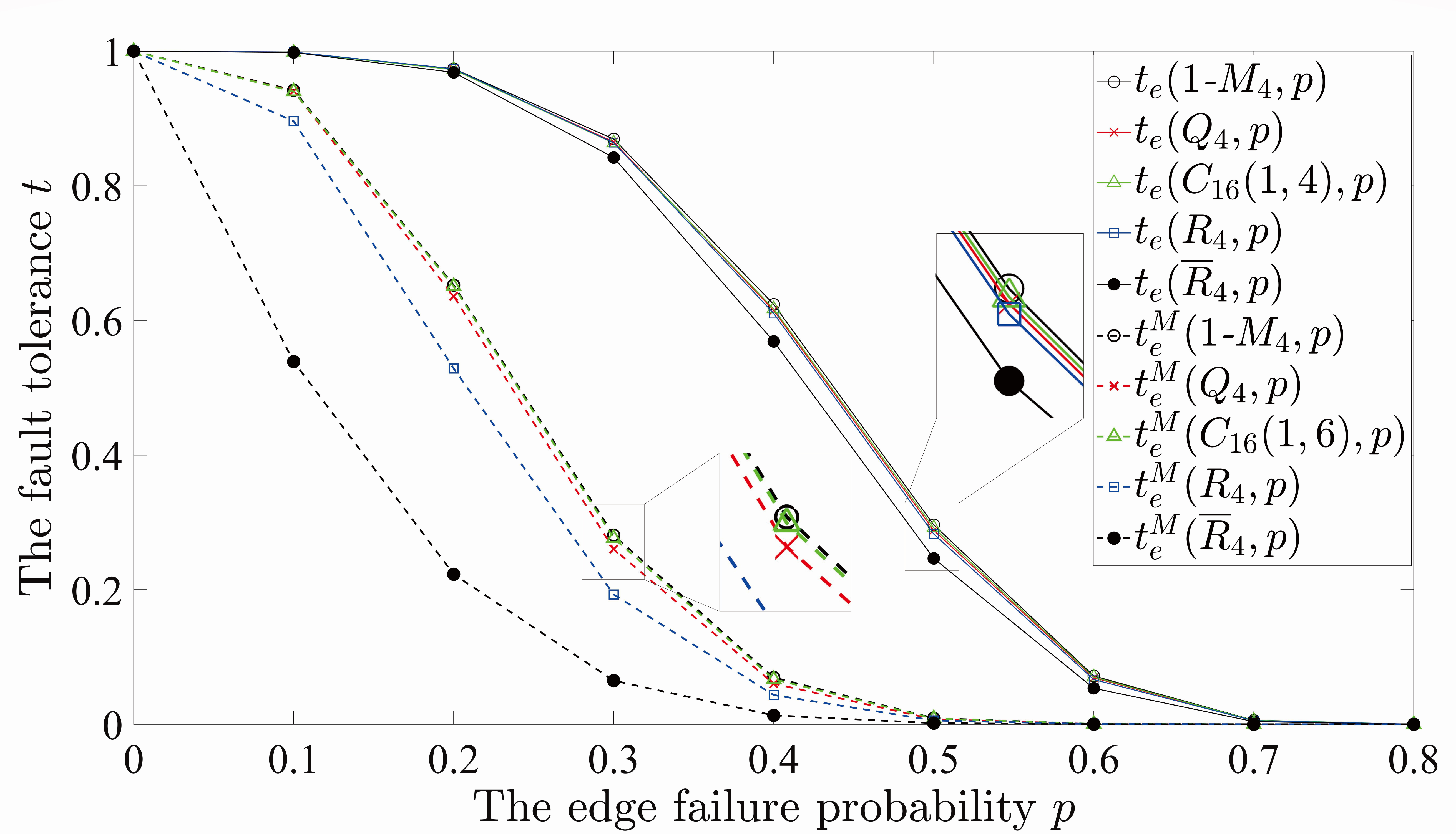}
\caption{The values of $t_e(G,p)$ (solid lines) and $t_e^M(G,p)$ (dotted lines) of some graphs with $p=0,0.1,\ldots,0.8$.}
\label{fig3}
\end{figure}

From TABLE \ref{tab1} and TABLE \ref{tab2} and Fig 3, we can observe that, compared with the 50 randomly generated graphs, the Hypercube-Like graphs $Q_4$, 0-$M_4$, 1-$M_4$ and the Circulant graphs  have both higher EF  tolerance and MEF tolerance. In particular, except the value of $t_e(G,0.1)$, 1-$M_4$ attains all the maximum values of both $t_e(G,p)$ and  $t^M_e(G,p)$ for every $p\in\{0.1,0.2,\ldots,0.9\}$ among all involved graphs. Even so,  it is a little surprising  that the differences of $t_e(G,p)$ among all the involved graphs (including the 50 randomly generated graphs) are very small.  In contrast, the differences of $t^M_e(G,p)$ are much larger.
%{\color{red}From Figure 3 we can also observe that a graph that has higher EF tolerance may also preserve higher MEF tolerance.}

We note that the upper bound on $t_e(G,p)$ given in Theorem \ref{lp} is very loose for not much large graphs. For an example, consider the probability that $Q_4$ is not connected when $p=0.1$. The upper bound gives $1-0.9992=0.0008$ while the numerical result gives $1-t_e(Q_4,0.1)=0.00162$, the former of which is about a half of the latter. Even so, this bound reveals  an asymptotical behavior for large scale regular graphs.  Recall that  if $p>1/2$ then  $t_e(Q_n,p)\rightarrow 0$, as $n\rightarrow+\infty$ (see Section III). On the other hand, our numerical result shows that the Hypercube $Q_4$ has higher EF tolerance than that of randomly generated graphs in general (see Figure 3 for an intuitive observation). This hints that if $p>1/2$ then  $t_e(G,p)\rightarrow 0$, as $n\rightarrow+\infty$, for almost all (in the sense of probability) $n$ regular graphs $G$ with $2^n$ vertices. This implies that, if the phase transition threshold for $G$ exists, then it should be not greater than $1/2$ for almost all $G$. Finally, we note that the graph attaining the maximum values of $t_e(G,p)$ and $t_e^M(G,p)$ among all involved graphs, i.e., the  M\"{o}bius Cube 1-$M_4$, is not symmetric (for the non-symmetry of  1-$M_4$, we refer to \cite{Cull}).

\section{Conclusion}
\noindent In this paper we considered the edge-fault tolerance of regular graphs in the sense of edge connectivity. In Section III we derived an upper bound on the EF tolerance $t_e(G,p)$ for regular graphs $G$, which reveals an asymptotical behavior when the scale of a graph is sufficiently large and $p$ is larger than certain values. In Section IV, we performed a simulation experiment on EF and MEF tolerance for some 4-regular graphs with 16 vertices. The numerical results show that, in addition to their well-structured properties for network models,  the Hypercubes, M\"{o}bius Cubes  and Circulant graphs, as reasonably expected, have also higher EF and  MEF tolerance in general. In particular, the M\"{o}bius Cube  1-$M_4$ has both the highest EF tolerance and  MEF tolerance among all involved graphs. Even so,  it is a little surprising  that the differences of the EF tolerances among all the involved graphs are very small. In contrast, the thing changes for MEF tolerance. This hints that, in contrast to  MEF tolerance,  the EF tolerance of regular graphs would be not strongly effected by the graph structure. Finally, we also point out that the graph attaining the maximum values of $t_e(G,p)$ and $t_e^M(G,p)$ among all involved graphs, i.e., the  M\"{o}bius Cube 1-$M_4$, is not symmetric. This implies that the symmetry of a graph is not a necessary structural feature that guarantees high fault tolerance of the graph.

To end the paper, we propose the following questions:
\vspace{2mm}\\
{\bf Question 2}. What is the ratio $\rho=t_e^M(Q_n,p)/t_e(Q_n,p)$ as $n\rightarrow +\infty$? In particular, is $\rho$ variant of the value of $p$?
\vspace{2mm}\\
{\bf Question 3}.  For any positive integer $n$ with $n\geq 4$, is it true that 0-$M_n$ or 1-$M_n$ attains  the maximum value of $t_e(G,p)$ (or  $t_e^M(G,p)$) among all Hypercube-Like graphs,  Ary Cubes and Circulant graphs of the same scale?

\section*{Declaration of competing interest}
 The authors declare that they have no known competing financial interests or personal relationships that could have appeared to influence the work reported in this paper.
 \section*{Acknowledgements}
This work was supported by the National Natural Science Foundation of China [Grant number, 12361070].

\end{document}